\documentclass [12pt,a4paper]{amsart}

\usepackage{amssymb,amsxtra,amsfonts}
\usepackage{epsf}              %
\usepackage{epsfig}             %

\openup\jot 
\newcommand\beq{\begin{equation}}
\newcommand\eeq{\end{equation}}

\newcommand{\IP}{\mathbb{P}}

\newcommand{\IC}{\mathbb{C}}
\newcommand{\IZ}{\mathbb{Z}}

\newcommand{\IT}{\mathbb{T}}

\newcommand{\M}{\mathcal{M}}

\newcommand{\cC}{\mathcal{C}}
\newcommand{\cL}{\mathcal{L}}

\newcommand{\Eps}{\mathcal{E}}

\newcommand{\cF}{\mathcal{F}}

\newcommand{\cO}{\mathcal{O}}

\newcommand{\cQ}{\mathcal{Q}}

\newcommand{\gl}{       \mathfrak{gl}     } %

\newcommand{\id}{\text{Id.}} 
\newcommand{\PVI}{$\text{P}_{\text{VI}}$}   %
\newcommand{\PV}{$\text{P}_{\text{V}}$}   %
\newcommand{\PIV}{$\text{P}_{\text{IV}}$}   %
\newcommand{\PIII}{$\text{P}_{\text{III}}$}   %
\newcommand{\Ptwo}{$\text{P}_{\text{II}}$}   %
\newcommand{\PI}{$\text{P}_{\text{I}}$}   %

\newcommand{\h}{\mathfrak{h}}

\newcommand{\Lia}{{\text{\rm Lie}}}     %

\newcommand{\pf}{\begin{bpf}}

\newcommand{\pfms}{\begin{bpfms}}
\newcommand{\epf}{\end{bpf}\hfill$\square$\\}           %
\newcommand{\epfms}{\end{bpfms}\hfill$\square$\\}               %

\newcommand{\idea}{\begin{bidea}}

\newcommand{\eidea}{\end{bidea}\hfill$\square$\\}           %

\newcommand{\sk}{\begin{bsk}}    %

\newcommand{\esk}{\end{bsk}\hfill$\square$\\}           %
\newcommand{\sketch}{\begin{bsketch}}%

\newcommand{\esketch}{\end{bsketch}\hfill$\square$\\}

\newcommand{\lbra}{{(\!(}}
\newcommand{\rbra}{{)\!)}}

\newcommand{\wt}{\widetilde}

\newcommand{\wh}{\widehat}
\newcommand{\al}{\alpha}

\newcommand{\be}{\beta}

\newcommand{\Ga}{\Gamma}
\newcommand{\la}{\lambda}
\newcommand{\La}{\Lambda}

\newcommand{\rank}{\text{\rm rank}}

\newcommand{\pr}{\text{\rm pr}}         %
\newcommand{\per}{\text{\rm perm}}         %

\newcommand{\sym}{\text{\rm Sym}}
\newcommand{\tr}{\text{\rm Tr}}

\newcommand{\Hom}{\text{\rm Hom}}
\newcommand{\Aut}{\text{\rm Aut}}

\newcommand{\SL}{\text{\rm SL}}

\newcommand{\GL}{\text{\rm GL}}
\newcommand{\PGL}{\text{\rm PGL}}

\newcommand{\Pic}{\text{\rm Pic}}

\newcommand{\End}{\text{\rm End}}
\newcommand{\diag}{{\text{\rm diag}}}

\newcommand {\eps}{\varepsilon}

\newcommand{\spq}{/\!\!/}

\def\mapright#1{\smash{
        \mathop{\longrightarrow}\limits^{#1}}}

\def\underset#1#2{\ \smash{\mathop{ #2 }\limits_{#1}}\ }

\theoremstyle{plain}
\newtheorem {hypo}{\bf\hspace{-\parindent}Hypothesis}

\newtheorem {thm}{Theorem}
\newtheorem {prop}[hypo]{Proposition}%

\newtheorem {lem}[hypo]{Lemma}%

\theoremstyle{definition}
\newtheorem {eg}[hypo]{Example}%
\newtheorem {exercise}{Exercise}%

\theoremstyle{remark}
\newtheorem {rmk}[hypo]{Remark}%

\begin{document}

\title[Quivers and difference Painlev\'e equations]{Quivers and difference Painlev\'e equations}
\author{Philip Boalch}
\address{\'Ecole Normale Sup\'erieure et CNRS\\
45 rue d'Ulm\\
75005 Paris\\
France} 
\email{boalch@dma.ens.fr}
\urladdr{www.dma.ens.fr/$\sim$boalch}

\begin{abstract}
We will describe natural Lax pairs for the 
difference Painlev\'e equations with affine Weyl
symmetry groups of types $E_6, E_7$ and $E_8$, showing that they do indeed arise as symmetries of certain Fuchsian systems of differential equations.
\end{abstract}

\dedicatory{
Dedicated to John McKay}

\maketitle

\renewcommand{\baselinestretch}{1.02}            %
\normalsize

\begin{section}{Introduction}
The theory of the Painlev\'e differential equations is now well established
and their solutions, the Painlev\'e transcendents, appear in numerous
nonlinear problems of mathematics and physics, from Einstein metrics to Frobenius manifolds.

However in certain applications these differential equations arise from a
limiting process, starting with a {\em discrete} nonlinear equation.
For example (see \cite{Sakai07}) the difference equation
$$f_{n+1}+f_n+f_{n-1}=\frac{an+b}{f_n} + c$$
arose in work of Br\'ezin and Kazakov \cite{BrezinKazakov} and Gross and Migdal \cite{GrossMigdal90long} in the calculation of a certain partition function of a $2$D quantum gravity model, and it becomes a differential Painlev\'e equation in a certain limit.
The classification and study 
of such ``discrete Painlev\'e equations'' 
has been an active area in recent years (see for example the reviews \cite{RamGra-review04} or \cite{Sakai07}).
On one hand many such equations have been found by studying a discrete analogue of
the Painlev\'e property, such as singularity confinement.
On the other hand 
Sakai \cite{Sakai-CMP01} has discovered a uniform framework to classify all second
order discrete Painlev\'e equations in terms of rational surfaces obtained by
blowing up $\IP^2$ in nine, possibly infinitely near, points.
In this framework, to each equation is associated an affine root system $R$
embedded in the affine $E_8$ root system, and thus an affine Dynkin diagram. (Upon removing the extending node, the Dynkin diagrams which occur correspond to the irreducible reflection subgroups of the $E_8$ Weyl group.) The (extended) affine Weyl group of the orthogonal 
$R^\perp$ of $R$ is then the symmetry group of the corresponding equation. 
The possible affine Weyl symmetry groups are roughly as follows:
\begin{figure}[h]
\epsfxsize=5.0cm
\input{sakai.sym.short2.pstex_t}
\end{figure}

\noindent
(See \cite{Sakai-CMP01} p.182 for the precise version.)
There are three types of discrete Painlev\'e equations:
difference, $q$-difference or elliptic and we will be concerned here solely with the (ordinary) difference equations---this rules out  the cases of $D_5, A_4$ and $A_2+A_1$.

It turns out that many of these difference Painlev\'e equations may be understood  as discrete symmetries of 
Painlev\'e differential equations: First the differential 
Painlev\'e equation
may be obtained as the equation controlling the isomonodromic deformations of a rational system of linear differential equations on $\IP^1$ of the form

\beq\label{eq: rat system}
\frac{d}{dz}- A
\eeq
for some matrix $A$ of rational functions.
Then the difference Painlev\'e equation arises as Schlesinger
transformations (or contiguity relations) 
of \eqref{eq: rat system}, obtained by applying rational gauge
transformation to shift the exponents at the poles by integers
(see \cite{Sakai-CMP01} pp.204-205).
This only applies to the bottom two rows of Sakai's list and the 
corresponding differential Painlev\'e equations are:

\begin{figure}[h]
\epsfxsize=5.0cm
\input{painleve.eqs.short2.pstex_t}
\end{figure}

\noindent
so for example the sixth Painlev\'e equation \PVI\  has symmetry group affine $D_4$ (which may be extended to affine $F_4$ by including diagram automorphisms).

Thus, looking at Sakai's table, this misses the cases with the richest symmetry groups,
namely we still have the cases of $E_6, E_7$ and $E_8$.

A basic question is thus to try to find linear systems \eqref{eq: rat system} with the desired symmetry groups, and living in two-dimensional moduli spaces.
(Explicitly this is Problem A of \cite{Sakai07}.)
A key observation is that to some extent it is irrelevant that there are corresponding differential 
Painlev\'e equations---as we will see the solution to this question turns out to be amongst certain Fuchsian systems which do not admit nontrivial isomonodromic deformations (in spite of the fact they are nonrigid).
Another point is that systems of higher rank appear. (Recall \cite{JM81} that all the differential Painlev\'e equations may be obtained by considering isomonodromy deformations of {\em rank two} systems \eqref{eq: rat system}.) 
Indeed the ranks will be given by the dimension of the largest irreducible representation of the finite group associated to the given Dynkin diagram via the McKay correspondence, viz:
$$\text{ rank $3,4,$ or $6$ for the cases of $E_6,E_7,$ or $E_8$}$$ respectively.

Said differently finding a linear system 
\eqref{eq: rat system} with the desired affine Weyl symmetry group 
(and living in a two-dimensional moduli space) amounts to finding a `Lax pair' for the corresponding  Painlev\'e difference equation (see \cite{RamGra-review04} pp.266 and 274; the Schlesinger transformations themselves give the other half of the Lax pair).
The fact that the system lives in a complex two-dimensional moduli space implies that, if written out explicitly in coordinates, the translation subgroup of the symmetry group constitutes a second order nonlinear difference equation.

Thus the aim of this article is to describe natural linear systems \eqref{eq: rat system} with the appropriate symmetries 
appearing in Lax pairs for the outstanding cases
of affine Weyl groups of types
$E_6, E_7$ and $E_8$ in Sakai's table.
The $E_8$ case is the  
most symmetric difference Painlev\'e equation.
(The $E_6$ case we will present here appeared recently from a different `Mellin transformed' viewpoint in Remark 2.14 of \cite{ArBo-duke}, which we are thus  generalising.) 
The origin of the desired spaces of Fuchsian systems lies in Kronheimer's work \cite{Kron.ale} which is very closely related to the McKay correspondence.
The main work to be done is to relate Kronheimer's work to Sakai's.

\vspace{3mm}

The layout of this article is as follows. 
Section 2 recalls the geometry of the Painlev\'e VI phase spaces as a motivating example. Section 3 then recalls the complex symplectic approach to Kronheimer's  asymptotically locally Euclidean hyperK\"ahler four manifolds, explaining how they arise as quiver varieties.
Section 4 recalls how quiver varieties for star-shaped quivers parameterise 
logarithmic connections on the trivial vector bundle over $\IP^1$ (i.e. Fuchsian systems). 
Then section 5 describes the spaces of Fuchsian systems corresponding to quivers whose underlying graph is a star-shaped affine Dynkin diagram. 
Then (families of) such spaces inherit a natural action of the corresponding finite Weyl groups. 
Section 6 then describes this action and how it extends to a birational action of the full affine Weyl group (the translation subgroup of which will constitute the desired difference Painlev\'e equation).
Section 7 ties things together by describing directly how to relate the spaces of Kronheimer to those of Sakai on which the nonlinear difference equation were originally defined (in essence one glues an affine line onto Kronheimer's spaces such that the affine Weyl group action becomes biregular rather than birational). 
Finally in the appendix we give a direct geometrical approach to the Weyl group symmetries (as well as some generalisations), addressing the philosophical question of ``where do the symmetries come from?" in terms of operations on connections.

The possibility of finding the $q$-difference analogue (or the elliptic analogue) of the Lax pairs of this
article will be studied elsewhere. 

\end{section}

\begin{section}{The Painlev\'e VI phase space}\label{sn: p6}

We will begin by describing the phase space of the sixth Painlev\'e
differential equation, in a slightly nonstandard way. (The corresponding contiguity relations lead to the
so-called d\PV\  equation which, as for the \PVI\ differential equation
itself, has symmetry group the affine Weyl group of type $D_4$.)

Painlev\'e VI arises as the equation controlling the isomonodromic deformations of $2\times 2$ Fuchsian systems of the form
\beq\label{eq: pvi system}
\frac{d}{dz}-A;\qquad A=\frac{A_1}{z} + \frac{A_2}{z-t} + \frac{A_3}{z-1}
\eeq
where $t\in\IC\setminus\{0,1\}$ is the deformation parameter.
Tensoring \eqref{eq: pvi system} by a logarithmic connection on a line bundle allows us to assume without loss of generality that each of the residues $A_i$, $i=1,2,3$ is a rank one matrix (i.e. we may add a scalar matrix to each $A_i$). By convention we choose $\theta_i\in\IC$ so that $A_i$ has eigenvalues $0,\theta_i$, $i=1,2,3$.
Finally we choose a fourth rank one matrix $A_4$ with eigenvalues $0,\theta_4$  such that $\sum_1^4 A_i$ is a scalar matrix:
\beq\label{eq: constraint}
A_1+A_2+A_3+A_4= \nu.1
\eeq
with $\nu\in\IC$, so that, upon taking the trace we see
\beq\label{eq: scalar constraint}
\sum_1^4\theta_i = 2\nu.
\eeq
(Note that the residue at infinity of $Adz$ is $A_4-\nu$ so in all cases 
$\theta_i$ is the difference, with appropriate sign, between the eigenvalues of the corresponding residue.)
Let $\cO_i\subset \gl_2(\IC)$ be the adjoint orbit of rank one matrices with eigenvalues $0,\theta_i$.
Each orbit may be identified with a coadjoint orbit using the trace pairing and so has a natural invariant complex symplectic structure, and the phase space we are interested in is the subspace of the product of the four orbits satisfying the constraint \eqref{eq: constraint}, modulo diagonal conjugation. But this is just the symplectic quotient of the product of orbits by diagonal conjugation at the value $\nu$ of the moment map, since the moment map is just the sum:
\beq\label{eq: phase space}
\cO_1\times\cdots\times\cO_4\underset{\nu}{\spq} \GL_2(\IC) = 
\{ (A_1,\ldots,A_4) \ \bigl\vert\ A_i\in\cO_i,\ \eqref{eq: constraint} \text{ holds }\}/\GL_2(\IC).
\eeq

Here the double slash ($\spq$) denotes the complex symplectic quotient and the single slash denotes the affine GIT quotient, taking the variety associated to the ring of invariant functions. (If $\theta_i=0$ we should take 
$\cO_i$ to be {\em closure} of the orbit of nonzero nilpotent matrices.)
These phase spaces are complex surfaces in general. Since any rank one matrix can be written as the tensor product of a vector and a linear form, we see that each orbit $\cO_i$ is itself the symplectic quotient of a vector space, the cotangent bundle of $\IC^2$, by the natural scaling action of $\IC^*$, at the value $\theta_i$ of the moment map:

$$\cO_i= T^*\IC^2\underset{\theta_i}{\spq}\IC^* = 
\{ (v,\al)\in \IC^2\times(\IC^2)^*\ \bigl\vert\ \al(v)=\theta_i \}/\IC^*.
$$ 

The map to $\cO_i$ is simply $(v,\al)\mapsto v\otimes \al$.
Combining these four $\IC^*$ quotients together with the above $\GL_2(\IC)$ quotient yields a description of the phase space \eqref{eq: phase space} as a symplectic quotient of a vector space $(T^*\IC^2)^4$ by the group $(\IC^*)^4\times \GL_2(\IC)$.

Now we know {\em a posteriori} that the set of these phase spaces admits a birational action of the affine $D_4$ Weyl group (due mainly to work of Okamoto \cite{OkaPVI}), with the set of parameters $\theta$  corresponding to the complex span of the $D_4$ root lattice. However the description above leads to a direct path from the $D_4$ root system to the phase spaces, as we will recall in the following sections. This will then generalise immediately to the root systems of types $E_6, E_7$ and $E_8$.

One can see an immediate link to $D_4$ since the condition for the phase spaces to be nonsingular is given by $\theta_i\neq 0$ (so there are no nilpotent orbits appearing) and $\sum_1^4\pm\theta_i\neq 0$ for any choice of signs  (this implies the system cannot be conjugated into triangular form
and thus that the action of $\PGL_2(\IC)$ is free). Now the $24$ vectors 
$\pm\epsilon_i, \sum_1^4\pm\epsilon_i/2$ make up a $D_4$ root system so the nonsingularity  conditions are just that the vector $\theta=\sum\theta_i\epsilon_i$ should not lie on any of the reflection hyperplanes for the $D_4$ Weyl group. The key facts (see \cite{Kron.ale}) to see the generalisation  are that for $\theta=0$ the space \eqref{eq: phase space} is the quotient $\IC^2/\Gamma$ (where $\Gamma\subset \SL_2(\IC)$ is the finite group corresponding to $D_4$ under the McKay correspondence), that the other spaces make up a semiuniversal deformation of this singularity, 
and that they have nice descriptions as symplectic quotients.

\begin{rmk}
Extrinsically the Okamoto transformations of \PVI\ were explained in \cite{AL-p6ims} in terms of the algebraic geometry of the moduli spaces, and intrinsically they were explained in \cite{k2p,srops} in terms of natural operations on connections, the Fourier--Laplace transformation and tensoring by connections on line bundles (which turned out \cite{Dett-Reit-p6} to be equivalent to the complex analytic version of Katz's middle convolution operation \cite{Katz-rls})---the appendix to the present article will extend this viewpoint.
\end{rmk}

\end{section}

\begin{section}{Kronheimer's ALE quiver varieties}  \label{sect, ALE}

We will recall the complex symplectic approach to Kronheimer's ALE spaces 
\cite{Kron.ale}, which are non-compact complex surfaces admitting 
 asymptotically locally Euclidean hyperK\"ahler metrics and which arise as deformations of Kleinian singularities. Subsequently these varieties were generalised 
by Nakajima \cite{Nakaj-quiver.duke}. (The article \cite{CB.ICM} of Crawley-Boevey has been useful to us here.)

Let $\cQ$ be a quiver (i.e. an oriented graph) whose underlying graph is a simply laced affine Dynkin diagram corresponding to a (finite) root system $R$ of rank $r$. 
Thus $R$ is of type $A_r$ with $r=1,2,3,...$, of type $D_r$ with 
$r=4,5,6,...$ or of type $E_r$ with $r=6,7,8$.
Let $I$ be the set of nodes and let $\cQ$ denote the set of edges.
The set $I$ corresponds to the simple roots $\eps_1,\ldots,\eps_r$ of $R$ plus $\eps_0$ corresponding to the extending node of the diagram.

Following McKay we associate a finite dimensional complex vector space $W_i$ to each node. Namely McKay observed that $\cQ$ may be identified with the 
McKay graph of a finite subgroup $\Ga$ of
$\SL_2(\IC)$ defined as follows. The nodes of the McKay graph correspond to the irreducible representations $V_i$ of $\Gamma$ (for each node $i\in I$) with $V_0$ being the trivial one-dimensional representation. Then nodes $i$ and $j$ are connected by $a_{ij}$ edges where $a_{ij}$ is the multiplicity of $V_j$ appearing in the tensor product of $V_i$ with the two dimensional defining representation $U$ of $\Gamma$:
$$U\otimes V_i \cong \bigoplus_j a_{ij} V_j.$$
(Since $U$ is self-dual $a_{ij}=a_{ji}$.)
Let $n_i:=\dim(V_i)$ be the dimension of the representation $V_i$ and let 
$\delta\in \IZ^I$ denote the vector of these dimensions.
Thus we may  associate a vector space $W_i:=\IC^{n_i}$ to each node. 
(Intrinsically $W_i$ should be thought of as the multiplicity space for $V_i$ in the regular representation of $\Gamma$---cf. \cite{Kron.ale}.)
In Lie theoretic terms the dimensions $n_i$ also arise as the coefficients expressing the longest root of $R$ in terms of the simple roots. That is the longest root
$$\rho= \sum_1^r n_i \eps_i$$
and $n_0=1$.
These coefficients $n_i$ are listed in \cite{BbkLie} or \cite{Kac-book} p. 54. For example for the root system of type $D_4$ the dimensions
are as in the figure below.

\begin{figure}[h]
\epsfxsize=5.0cm
\input{dynkinD4.pstex_t}

\text{Extended $D_4$ Dynkin diagram}

\text{ }\vspace{-4mm}

\text{McKay graph of the quaternion group $\{\pm1,\pm i,\pm j,\pm k\}$}
\end{figure}

Now to construct the quiver variety we consider  for each edge of $\cQ$ a homomorphism in each direction between the vector spaces at each end:

$$W=W(\cQ):=\bigoplus_{e\in Q} \Hom(W_{t(e)},W_{h(e)}) \oplus  \Hom(W_{h(e)},W_{t(e)})$$
where, for a directed edge $e\in \cQ$, $h(e),t(e)\in I$ are its head and tail respectively.

The vector space $W$ inherits a linear complex symplectic structure if we identify  $\Hom(W_{t(e)},W_{h(e)}) \oplus  \Hom(W_{h(e)},W_{t(e)})$ with the cotangent bundle  $T^*\Hom(W_{t(e)},W_{h(e)})$ in the obvious way.

Now consider the group
$$\wt G=\wt G(\cQ):= \prod_{i=0}^r\GL(W_i)$$
which acts on the sum of the vector spaces $W_i$ and thus has an induced action on $W$. (If $\phi\in \Hom(W_i,W_j)$ and $g=(g_0,\ldots,g_r)\in \wt G$ 
then $g\cdot \phi = g_j \circ \phi \circ g_i^{-1}$.)

The group $\wt G$ has a central subgroup $\IT$ isomorphic to $\IC^*$ which acts trivially on $W$ (simply take each $g_i$ to be the same nonzero scalar).
Define $G = \wt G / \IT$ to be the quotient group (which acts effectively on $W$). Thus $G$ has a centre $Z$ of dimension $r$ (obtained by setting each $g_i$ to be an arbitrary nonzero scalar).

A moment map $\mu = (\mu_0,\ldots,\mu_r): W \to \prod_{i=0}^r\gl(W_i)$ for the action of $\wt G$ on $W$ is given by 
$$\mu_i(\{\phi\}) := \sum_{e\in Q, h(e)=i} \phi_e \circ \phi_{e^*} -
\sum_{e\in Q, t(e)=i} \phi_{e^*} \circ \phi_{e}$$
where $\phi_e\in \Hom(W_{t(e)},W_{h(e)})$ and $e^*$ is the edge $e$ with reversed orientation.

Taking the sum of the traces of the components of $\mu$ clearly yields zero (since $\tr(AB)=\tr(BA)$ for any compatible rectangular matrices $A,B$) and the subspace on which this sum of traces is zero is naturally identified with the (dual of the) Lie algebra of $G$, and so $\mu$ is also a moment map for $G$.

Finally the corresponding quiver varieties are obtained by performing the complex symplectic quotient of $W$ by $G$ at any central value $\lambda$ of the moment map $\mu$:
$$N_\cQ(\lambda) = N_\cQ(\lambda, \delta) 
= W\underset{\la}{\spq} G = \mu^{-1}(\lambda)/G$$
where $\delta$ is the dimension vector whose components were the dimensions $n_i$.
(We may also view this as $W\spq_\la\wt G$ since $\IT$ acts trivially.)

Let $\h$ be the set of such coadjoint invariant elements (i.e. central elements) $\lambda$ of the dual of the Lie algebra of $G$. Since the trace pairing relates the adjoint action and the coadjoint action, $\h$ may be identified with the Lie algebra of the central subgroup $Z$ of $G$ and so is a complex vector space of dimension $r$. Concretely it is given as a subspace of 
$\prod_{i=0}^r\gl(W_i)$ by the  $r+1$-tuples of scalar matrices whose sum of traces is zero. In other words if each diagonal entry of the $i$th matrix is 
$\lambda_i$ then the constraint is that
\beq\label{eq: de.la=0}
\delta\cdot\lambda := \sum_0^r n_i \lambda_i = \sum_0^r\tr_{W_i}(\la_i)= 0.
\eeq
Carrying this out for the case of $D_4$ should now be familiar from the previous section, and note in particular that the above constraint matches up with that of \eqref{eq: scalar constraint}.

\end{section}

\begin{section}{Quivers and matrices}\label{sn: quiver+matrices}

Of course one can carry out the construction of the previous section for any quiver and any dimension vector, although in general these spaces will probably not be relevant to Painlev\'e equations or isomonodromy in general.
However there are many cases which are: the  phase spaces of the Schlesinger equations for $\GL_n(\IC)$ arise as quiver varieties for {\em star-shaped} quivers.

This relationship is by now well-known and has been used in an essential way by Crawley-Boevey (see \cite{CB.ICM}) in his work on the Deligne--Simpson problem.

The $\GL_n(\IC)$ Schlesinger phase spaces are just the moduli spaces of logarithmic connections on trivial rank $n$ holomorphic vector bundles over the Riemann sphere, obtained by fixing the adjoint orbits of the residues (the name ``moduli spaces of Fuchsian systems" seems most suitable).
Choosing a trivialisation of the bundle and a coordinate $z$ on $\IP^1$, such a connection will be of the form 
\beq\label{dis: fuchsian systems}
d - Adz,\qquad A=\sum_{i=1}^{m-1}\frac{A_i}{z-a_i}
\eeq
which is equivalent to the Fuchsian systems $\frac{d}{dz}-A$.
If we set $a_m=\infty$ and take $A_m$ to be the residue of $Adz$ at $\infty$
then the residue theorem says that
$$\sum_1^m A_i =0.$$
If we fix the residues to be in chosen (co)adjoint orbits $A_i\in \cO_i$ then we notice that $\sum_1^m A_i$ is just the moment map for the diagonal conjugation action of $\GL_n(\IC)$ on the product of orbits, and that this action corresponds to the gauge equivalence of such connections, so the moduli space is just the symplectic quotient
\beq\label{display: fuchsian spq}
\cO_1\times \cdots \times \cO_m \spq \GL_n(\IC)
\eeq
at the zero value of the moment map (cf. \cite{Hit95long, smid}).
The rank two, four pole case is the Painlev\'e VI case described in a slightly different way in section \ref{sn: p6}, which we have now recognised as $D_4$ quiver varieties. However the more general case \eqref{display: fuchsian spq}
can also be obtained from a star-shaped quiver, as follows (see \cite{CB.ICM}).

Let $w_i$ be the degree of the minimal polynomial $p_i$ of $A_i\in \cO_i$ and let $\xi_{ij}$ be the roots of $p_i$ ($i=1,\ldots,m$ and $j=1,\ldots,w_i$).
Let $\cQ$ be the star-shaped quiver with a central vertex labelled $*$ and $m$ legs of lengths $w_i-1$ for $i=1,\ldots,m$. 
Label the $j$th vertex of the $i$th leg  by $ij$ with $1\le j < w_i$ (with $j$ increasing as we go away from $*$). Orient all the edges to point towards the central vertex $*$.
The corresponding dimension vector is obtained by assigning dimension $n$ to the central vertex and the dimension
$$n_{ij}:= \rank(A_i-\xi_{i1})(A_i-\xi_{i2})\cdots(A_i-\xi_{ij})$$
to the vertex $ij$.
This specifies the quiver and dimension vector corresponding to the given orbits.

\begin{exercise}\label{eg: legex}
Consider the following leg of a star-shaped quiver, with the given vector spaces at each node.

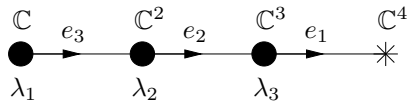
\begin{figure}[h]
\epsfxsize=5.0cm
\input{legexercise.pstex_t}

\end{figure}

i) Let $\la_i\in\End(\IC^i)$ be a scalar matrix and write down the 
condition for each component $\mu_i$ of the moment map $\mu$ to take the value $\la_i$ for $i=1,2,3$.

ii) If $\mu(\{\phi\})=\lambda$, deduce the eigenvalues of the matrices
$$\phi_{e^*_i}\circ\phi_{e_i},\qquad\phi_{e_i}\circ\phi_{e_i^*}$$
for $i=3,2,1$.

iii) Thus observe that $\phi_{e_1}\circ\phi_{e_1^*}\in\End(\IC^4)$
has eigenvalues $$0,\ -\la_3,\  -\la_3-\la_2,\  -\la_3-\la_2-\la_1$$
and so the $\la_i$'s arise as eigenvalue {\em differences}.
\end{exercise}

\begin{rmk}{(Eigenvalue conventions.)\ }
A set of parameters $\la$ for a star-shaped quiver satisfying $\la\cdot \delta=0$ thus specifies (and orders) the 
eigenvalues of all of the residues of the corresponding Fuchsian system 
\eqref{dis: fuchsian systems}, modulo tensoring by a scalar system of the form
$d-\sum_1^{m-1} \La_i\frac{dz}{z-a_i}$ with $\La_i\in\IC$. In other words, up to adding a scalar $\La_i$ to $A_i$ for $i=1,\ldots,m-1$ (and varying the residue $-\sum_1^{m-1}A_i$ at infinity in the corresponding way). There are two standard ways of normalising the residues to remove this indeterminacy, and both are useful:

(1) ``Trace zero": choose each $\La_i$ such that each residue has zero trace,

(2) ``Determinant zero": choose $\La_i$ such that the first eigenvalue of each residue $A_1,\ldots,A_{m-1}$ is zero---it is this normalisation that arises in the first instance from the quiver viewpoint. (The residue at infinity is then $A_m-\nu$ where $\nu$  is minus the component of $\lambda$ associated to the central node of $\cQ$,  $\sum_1^m A_i=\nu$, and the first eigenvalue of $A_m$ is also zero.) 

It is a triviality to go between these two normalisations, but one should be aware of them.   
\end{rmk}

\end{section}

\begin{section}{Star-shaped affine Dynkin diagrams}

Now we will look at the area of intersection of the two previous sections, i.e. the case of a star-shaped affine Dynkin diagram. This leaves just the cases of $D_4$ and $E_6, E_7, E_8$ from the original $ADE$ list. Inverting the construction of the previous section we can then view Kronheimer's $E$-type ALE spaces as moduli spaces of Fuchsian systems on $\IP^1$, as we have already done for $D_4$.

The $E$-type affine Dynkin diagrams and dimension vectors, which are the same as the Platonic McKay graphs, are as in the following figures.

\begin{figure}[h]
\epsfxsize=5.0cm
\input{dynkinE6.pstex_t}

\text{Extended $E_6$ Dynkin diagram}

\text{ }\vspace{-4mm}

\text{McKay graph of the binary tetrahedral group}
\end{figure}

\begin{figure}[h]
\epsfxsize=5.0cm
\input{dynkinE7.pstex_t}

\text{Extended $E_7$ Dynkin diagram}

\text{ }\vspace{-4mm}

\text{McKay graph of the binary octahedral group}
\end{figure}

\begin{figure}[h]
\epsfxsize=5.0cm
\input{dynkinE8.pstex_t}

\text{Extended $E_8$ Dynkin diagram}

\text{ }\vspace{-4mm}

\text{McKay graph of the binary icosahedral group}
\end{figure}

In all these cases there are just three legs so they correspond to Fuchsian systems with just three poles on $\IP^1$. Up to coordinate transformations we can always choose to put the poles at $0,1,\infty$ and so the systems are of the form
$$\frac{d}{dz} - \left(\frac{A_1}{z}+\frac{A_2}{z-1}\right)$$
and the residue at infinity is $A_3=-(A_1+A_2)$.
(Thus in general these systems have no nontrivial isomonodromic deformations, only discrete symmetries.) 
Now given any value of the parameters $\lambda$ the prescription of the previous  section applied to these graphs and dimension vectors determines the size and orbit type of the three matrices $A_i$ that occur. 
When the parameters are generic the matrices are semisimple with eigenvalues distributed as in the following table.

\begin{center}
\begin{tabular}{| c |c |c |c | c |}\hline
      &  Size $n$ & $\cO_1$ & $\cO_2$ & $\cO_3$ \\\hline
$E_6$ & $3$ &  $1,1,1$ & $1,1,1$ & $1,1,1$ \\
$E_7$ & $4$ &  $2,2$ & $1,1,1,1$ & $1,1,1,1$ \\
$E_8$ & $6$ &  $3,3$ & $2,2,2$ & $1,\ldots,1$\\\hline
\end{tabular}
\end{center}

The last three columns list various partitions of the size $n$ and indicate the eigenvalue multiplicities. Thus for example $2,2,2$ on the third row indicates orbits of $6\times 6$ matrices with three distinct eigenvalues each repeated twice. Thus the corresponding quiver varieties are isomorphic to complex symplectic quotients
$$\cO_1\times\cO_2\times\cO_3\spq \PGL_n(\IC)$$
with the orbits $\cO_i$ as indicated in the table.

In the $E_6$ case we recognise the generic three-poled $3\times 3$ systems mentioned in \cite{ArBo-duke} Remark 2.14 (it would be interesting to compute the Mellin transforms of the $E_7$ and $E_8$ systems in the table above to obtain ``difference-difference" Lax pairs).

\begin{exercise}
Confirm directly in each case that:

\nobreak
i) the dimension of the space of parameters (specifying the eigenvalues of the matrices) is equal to the rank of the finite root system,

ii) that the above symplectic quotient does indeed have expected complex dimension $2$, and 

iii) that as in the $D_4$ case above, if the parameters $\lambda$ are not on any of the root hyperplanes then the triple $A_1,A_2,A_3$ (with $A_i\in \cO_i$)
is irreducible (i.e. cannot be simultaneously conjugated into a block triangular form) and also there are no extra coincidences between the eigenvalues of each $A_i$.

{
\renewcommand{\baselinestretch}{0.85}              %
\Small
(For example for $E_8$: i) we need to choose $11$ eigenvalues but there is one relation since the sum of the traces should be zero, and two more parameters are irrelevant since we can tensor with a scalar system with poles at 
$0,1,\infty$. For ii) the product of orbits has dimension $72=18+24+30$ and $\PGL_6(\IC)$ has dimension $35$ so the expected dimension is 
$72-2\times 35 = 2$. For iii) write the three sets of eigenvalues as 
$$(a,a,a,b,b,b),\quad (c,c,d,d,e,e), \quad (f,g,h,i,j,k).$$
If the $A_i$ can be put in a block triangular form with an {\em irreducible} $r\times r$ block on the diagonal, then taking the sum of the traces of these blocks gives a linear relation amongst a subset of the eigenvalues:
For $r=1$ or $5$ there are $2\times 3 \times 6= 36$ possibilities,
for $r=2$ or $4$ there are $1\times { 3 \choose 2 }\times { 6 \choose 2}=45$ possibilities, and for $r=3$ there are $2\times 1\times {6\choose 3}=40$ possibilities\footnote{
Note for example that a $3\times 3$ block with eigenvalues of the form 
$(a,a,b), (c,c,d), (f,g,h)$ is necessarily reducible so not counted in the $r=3$ case.}, giving altogether $202$ linear relations (each appearing together with its minus so yielding $101$ hyperplanes in the parameter space).
If the parameters are off of these hyperplanes then clearly no such block triangular decomposition is possible.
Further, for the eigenvalues to not have extra coincidences yields $1+3+15=19$ further hyperplanes, or $38$ linear relations taking both possible signs of each relation. These $202+38=240$ linear relations are the $E_8$ roots.)

\renewcommand{\baselinestretch}{1}
\small
}
\renewcommand{\baselinestretch}{1}
\normalsize

\end{exercise}

\begin{rmk} This list has also appeared in work of Kostov \cite{Kostov-0idx} (which we learnt of from Simpson \cite{Simpson-katz}), but the link with root systems was not made.
\end{rmk}

\begin{rmk} Observe (as in \cite{Kron.ale}) that since the vector spaces attached to the nodes are of the same dimensions as the irreducible representations of the corresponding binary group, there is a simple formula for the complex dimension of the vector space $W$ from which the quiver variety was constructed: its complex dimension is twice the order of the binary group (e.g. $240$ for $E_8$).  
\end{rmk}

\begin{exercise} (Non-star-shaped examples and irregular singularities.)
Relate the symplectic quotient descriptions of the quiver varieties of affine Dynkin types 
$A_3,A_2$ and $A_1$ to the symplectic quotient descriptions (\cite{smid} section 2) of the moduli spaces of rank two (irregular) linear systems on the Riemann sphere, with poles of orders $(2,1,1), (3,1)$ and $(4)$ respectively. Hints:

i) Use an extended orbit (denoted $\wt \cO\cong T^*G\times \cO_B$ in \cite{smid}) at each irregular singularity. Thus after performing the $G:=\GL_2(\IC)$ symplectic quotient
the spaces are respectively the symplectic quotients of
$$ \cO_B^{(2)}\times\cO_1\times\cO_2,\qquad \cO_B^{(3)}\times\cO_1,\qquad  \cO_B^{(4)}$$
by the diagonal conjugation action of the torus
$T=\left\{\left(\begin{smallmatrix} * & \\ & 1 \end{smallmatrix}\right)\right\} \subset G$,
where the $\cO_i$ are coadjoint orbits of $G$, and the $\cO_B^{(k)}$
are coadjoint orbits of a nilpotent group (the group $B_k$ of based $k$-jets of bundle automorphisms).

ii) Show directly that $\cO_B^{(k)}\cong T^*\IC^{k-2}$ as Hamiltonian $T$-spaces, where $T\cong\IC^*$ acts on $\IC^{k-2}$ by scaling. (This is tricky for $k=4$, but possible.)

iii) Write out the corresponding ALE symplectic quotients explicitly in these cases, and  symplectic quotient by a  suitable codimension one subtorus such that the remaining $\IC^*$ quotient is precisely that of 1) above (up to tensoring to make the $\cO_i$ orbits of rank one matrices).
\end{exercise}

\end{section}

\begin{section}{Affine Weyl group action}

The next task is to show that the Fuchsian systems appearing in the previous section do indeed have the corresponding affine Weyl group symmetries. Fortunately the finite Weyl group action is now quite well-known (and was known to Kronheimer) and the translations come from Schlesinger transformations.

Let $\cQ$ be a star-shaped quiver (for example as considered in Section 
\ref{sn: quiver+matrices}). (This is a more general context than strictly necessary here but these more general symmetries, for example changing the dimension vector, will give many interesting identifications between various spaces of Fuchsian systems.)

One can then define a Weyl group as follows.
Let $I$ be the set of nodes and let $r+1=\#I$ be the number of nodes.
Define the $(r+1)\times (r+1)$ Cartan matrix to be 
$$C=2\ \id - A$$
where the $i,j$ entry of $A$ is the number of edges connecting the $i$th and  $j$th nodes. 
The {\em root lattice} $\IZ^I=\bigoplus_0^r \IZ\eps_i$ inherits a bilinear form defined by 
$$(\eps_i,\eps_j) = C_{ij}.$$

In turn we can define simple reflections $s_i$ acting on the root lattice by the 
formula
$$s_i(\be) := \be-(\be,\eps_i)\eps_i$$
for any $i\in I$.
The Weyl group of the quiver is then the group generated by these simple reflections.
We can also define dual reflections $r_i$ acting on the vector space 
$\IC^I$ by the formula 
$$r_i(\la) = \la - \la_i \al_i$$
where $\la = \sum_0^r \la_i\eps_i\in \IC^I$ with $\la_i\in\IC$ and where for each $i\in I$
$$\al_i := \sum_{j\in I} (\eps_i,\eps_j)\eps_j\in \IC^I.$$ 
By construction one has that $s_i(\be)\cdot r_i(\la) = \be \cdot \la$,
where the dot denotes the pairing given by 
$\eps_i\cdot\eps_j = \delta_{ij}$.
The main result we wish to quote is then
\begin{thm}[\cite{Nakaj-quiver.duke, Rump, CB-H, Nakaj-refl}]
Let $\cQ$ be a  fixed quiver as above.
Then, if $\la_i\ne 0$, 
there is a natural isomorphism between the the quiver variety with dimension vector $\be$ and parameters $\la=\sum\la_i\eps_i$ 
and that with dimension vector $s_i(\be)$ and parameters $r_i(\la)$:
$$N_{\cQ}(\la,\be) \cong N_{\cQ}(r_i(\la),s_i(\be)).$$
\end{thm}

The desired reflection functors are constructed in \cite{CB-H} section $5$ (see also \cite{Nakaj-quiver.duke, Nakaj-refl, Rump}) 
and that they induce isomorphisms of the moduli spaces is stated explicitly \cite{CB-decomp} Lemma 2.1 and proved in \cite{CB-mmap} Lemma 2.2.
(Apparently they were first written down independently by Nakajima 
\cite{Nakaj-japaneses.refls} and Rump \cite{Rump}.)\footnote{Of course for $D_4$ they were in some sense found earlier by Okamoto \cite{OkaPVI} (but even the most interesting Okamoto transformation was detected earlier at a quantum level by Regge \cite{Reg-6j}, cf. \cite{roks}).}

Now let us specialise to one of the four cases $D_4, E_6, E_7, E_8$ we are interested in. (This Weyl group action will then give all the desired symmetries except the translations.)
The generating reflections will be described geometrically in the appendix; here
we will check that the group acting is indeed the finite Weyl group of the corresponding type (i.e. $D_4, E_6, E_7$, or $E_8$).

Thus $\cQ$ is an extended Dynkin diagram and we label the nodes by $I=\{0,\ldots,r\}$ with 
$0$ labelling the extending node (i.e. any of the extremal nodes with dimension $1$). Let $\delta\in\IZ^I$ be the corresponding dimension vector (i.e. the dimensions written by the nodes in the figures above).

Then $\delta$  spans the kernel of the bilinear form $(\ ,\ )$ and so in particular all the reflections $s_i$ fix 
$\delta$. Recall from \eqref{eq: de.la=0} that the space $\h$ of coadjoint invariants is just $\h=\{ \la \ \bigl\vert\ \la\cdot\delta=0\}$.
Thus if $\la_i\ne 0$  the above construction yields an isomorphism
$N_{\cQ}(\la)\cong N_{\cQ}(r_i(\la))$. 
We will check that this generates the finite Weyl group, so that we have an action of the  Weyl group on the family of quiver varieties parameterised by the regular values of $\la$.

First note that each $\al_i:=  \sum_{j\in I} (\eps_i,\eps_j)\eps_j$ is actually in $\h\subset \IC^I$ and that $\al_1,\ldots,\al_r$ form a basis of $\h$.
Then define a bilinear form on $\h$ by declaring
$$\lbra\al_i,\al_j\rbra = C_{ij}\qquad\text{for $i,j=1,\ldots,r$}.$$
This is non-degenerate since $C$ (with indices $\ge 1$) is the Cartan matrix of a (non-affine) Dynkin diagram, and so $\al_1,\ldots,\al_r$ make up the simple roots of a copy of the corresponding finite root system.  
Now by definition the simple reflections of this root system are of the form
$$\la \mapsto \la - \lbra\la,\al_i\rbra\al_i.$$
But this is just the restriction of $r_i$ to $\h$ because:
\begin{lem}\label{lem: ips}
For all $\la\in \h$ we have $\lbra\la,\al_i\rbra = \la \cdot \eps_i =: \la_i$. 
\end{lem}
\pf
On basis vectors: $\lbra\al_j,\al_i\rbra=C_{ji}=\left(\sum_k C_{jk}\eps_k\right)\cdot \eps_i=\al_j\cdot\eps_i$. 
\epf

Hence the group of automorphisms of $\h$ induced by the Weyl group generated
by $r_1,\ldots,r_r$ is the desired finite Weyl group. (Note that $r_0$ gives nothing more since $\al_0$ is minus the longest root; the fact that $\delta$ is in the kernel of $(\ ,\ )$ implies $\sum_0^rn_i\al_i=0$.)

Next we will verify  that the translations arising from Schlesinger
transformations generate the translation subgroup of the corresponding
{\em extended} affine Weyl group. Recall that for a simply laced root system
$R$ (as we have here) the translation subgroup of the affine Weyl group is the
root lattice $Q(R)= \langle \al_1,\ldots,\al_r \rangle_\IZ\subset \h$ and the translation
subgroup of the extended affine Weyl group is the finer lattice
$$P(R) = \{ \la\in\h \ \bigl\vert\ \lbra \la,\al_i \rbra\in \IZ \text{ for all
$i=1,\ldots,r$}\}\supset Q(R).$$

Now by definition the Schlesinger transformations are the transformations of
the Fuchsian systems
induced by rational gauge transformations with poles at the poles of the
system (see \cite{JM81} for more details and some explicit formulae).
Thus they vary the eigenvalues of the systems by integers and induce
birational isomorphisms between the moduli spaces.
Of course we wish to obtain systems with eigenvalue distributions of the same
type, so the possible integer shifts are restricted. In the notation we are
using the Schlesinger transformations yield the integer shifts given by the
lattice
$$\{ \la \ \bigl\vert\ \la_i\in\IZ, \la\cdot\delta=0\}.$$
Clearly, by Lemma \ref{lem: ips}, this is the desired lattice $P(R)$.

Thus altogether we obtain the full extended affine Weyl group of birational symmetries of the corresponding Fuchsian systems. (In the $D_4$ and $E_6$ cases there are further symmetries corresponding to the automorphisms of the corresponding non-extended Dynkin diagram, which may be obtained by suitable automorphisms of the underlying Riemann sphere on which the Fuchsian system lives.)

\begin{rmk}
Another way of thinking of the Schlesinger transformations, at least for generic parameters (no nonzero integer eigenvalue differences), is as follows.
Consider the moduli space $\M(\la)$ of logarithmic connections on degree zero vector bundles over $\IP^1$, with poles and residue types as for the Fuchsian systems above. Thus the Fuchsian system moduli space $\M^*(\la)\cong \cO_1\times\cO_2\times\cO_3\spq \PGL_n(\IC)$  is a subset of $\M(\la)$ (and presumably arguments with determinant line bundles will show it is the complement of a divisor).
Then the Riemann--Hilbert map furnishes an analytic isomorphism 
$$\M(\la) \mapright{\cong} M(\la):=
\{ (M_1,M_2,M_3)\in\cC_1\times \cC_2\times\cC_3\ \bigl\vert\ 
M_1M_2M_3=1 \}/\PGL_n(\IC)$$
with the corresponding space of monodromy data. (Here $\cC_j\subset \GL_n(\IC)$ is the conjugacy class $\exp(2\pi i\cO_j)$.)
But it is immediate that if $\nu\in P(R)$ then $M(\la)=M(\la+\nu)$ and so using the Riemann--Hilbert map we obtain analytic isomorphisms 
$$\M(\la)\cong \M(\la+\nu).$$
But we know this map is algebraic on the open subset $\M^*(\la)$ and so is a (biregular) algebraic isomorphism.
This map extends the Schlesinger transformations since they 
clearly preserve the monodromy data.
(These isomorphisms may also be described algebraically in terms of elementary transformations, which are described for example in \cite{EV-98}.)
\end{rmk}

\end{section}

\begin{section}{Blow-ups of the projective plane in $6,7,8$ or $9$ points}

Finally, to tie things together, we will describe how the $E_6, E_7$ and $E_8$ ALE spaces of Kronheimer described above are related to the corresponding spaces of Sakai \cite{Sakai-CMP01}, on which the difference equations were originally defined.
The key point is that one is essentially adding an affine line (and the biregular Weyl group action then extends to a biregular affine Weyl group action).

First recall, in the case of $E_r$, that  Kronheimer's construction above yields a family of affine surfaces parameterised by $\IC^r$, for $r=6,7,8$, and that these surfaces will be non-singular if the parameters do not lie on one of the reflection hyperplanes of the (finite) Weyl group action on 
$\IC^r$. 
Kronheimer shows (\cite{Kron.ale} section 4) that the quotient of this family of varieties (denoted $Y\to Z\otimes \IC$ in \cite{Kron.ale}) by the Weyl group is the semi-universal deformation of the corresponding Kleinian singularity (with $\IC^*$ action).
In particular the nonsingular fibres are isomorphic to the corresponding fibres of the simultaneous resolution of this deformation
constructed by Brieskorn, Slodowy and Tjurina. Tjurina's approach to these fibres \cite{tjurina} was described by Pinkham (\cite{Pinkham} pp. 190-197): they are obtained by blowing up $r$ points in the smooth locus of a cuspidal cubic in $\IP^2$ (and then removing the strict transform of the cubic). Since this smooth locus is isomorphic to $\IC$ the parameter space is again a copy of $\IC^r$. The picture to have in mind is that if all $r$ points are equal to the inflection point of the cubic (where the tangent line has third order contact) then there will be a tree of $r$ $(-2)$-curves in the blow-up, intersecting according to the $E_r$ Dynkin diagram, since the tangent line becomes a $(-2)$-curve on the third blow-up.
(This tree then becomes the $E_r$ singularity upon taking the anti-canonical model of the surface, and deforming the $r$ points gives the semi-universal deformation.)

On the other hand the corresponding spaces of Sakai are obtained as follows:\footnote{The cases $E_8,E_7,E_6$  are labelled by ``Add. 1-3" or by $A_0^{(1)**}, A_1^{(1)*}, A_2^{(1)*}$ respectively in \cite{Sakai-CMP01}.  }

\noindent$\bullet$\  $E_8$-case: blow-up $9$ points on the smooth locus of a cuspidal cubic in $\IP^2$,

\noindent$\bullet$\  $E_7$-case: blow-up $9$ points on the smooth locus of the degenerate cubic in $\IP^2$ formed by a line {\em touching} a conic ($3$ points on the line and $6$ on the conic),

\noindent$\bullet$\  $E_6$-case: blow-up $9$ points on the smooth locus of the degenerate cubic in $\IP^2$ formed by three lines meeting in a single point ($3$ points on each line).

In each case one should finally remove the strict transform of the cubic, and there is a genericity condition (that there should be no other cubic through the same $9$ points). 
Quotienting by the (identity component) of the subgroup of $\Aut(\IP^2)=\PGL_3(\IC)$ preserving these cubics again identifies the spaces of parameters with $\IC^r$ in each case (see \cite{Sakai-CMP01} for explicit coordinates). 
In each case the resulting family of surfaces admits a biregular action of the corresponding (extended) affine Weyl group (at least over the complement of the reflection hyperplanes), described explicitly in terms of Cremona transformations---in other words the action just changes the way a given surface is expressed as a blow-up of $\IP^2$. 

Now in the $E_8$ case it should be clear how to relate the surfaces of Sakai and Tjurina: one is either blowing up $8$ or $9$ points in the smooth locus of a cuspidal cubic and then removing the strict transform of the cubic. Thus, in general, the difference is the affine line given by the exceptional curve of the ninth point minus its intersection point with the strict transform of the cubic.
(In this $E_8$ case the genericity condition is equivalent to the $9$ points {\em not} summing to zero in the group law of the cubic and up to automorphisms one may fix the value of the sum of the $9$ points to be any arbitrary nonzero value, and label it $1$. Thus up to automorphisms we see that the space of parameters is just the choice of the first $8$ points.
Similar statements hold in the other cases---see \cite{Sakai-CMP01}.)

In the $E_7$ and $E_6$ cases it is trickier to relate the surfaces of Sakai and Tjurina, since we are blowing up different numbers of points on different cubics. However in the end the difference still amounts to an affine line in the exceptional divisor 
of the blow-up of one extra point, as follows. 

A useful fact we will recall first is that if  a cuspidal cubic is blown-up at the cusp then one obtains two rational curves touching at a point. If this point is further blown-up then one obtains three rational curves intersecting at a single point (cf. \cite{Hartshorne} 3.9.1).

Now for the $E_7$ case the key point is to see how 7 points on the smooth locus of a cuspidal cubic lead to 6 points on a conic in $\IP^2$ and two points on a line 
that touches it. (Then to get to Sakai's picture we need only choose one further point on the line, and modulo automorphisms there is actually no choice involved here.) 
To do this one does the standard quadratic Cremona transform based at the cusp $c$ and at two of the seven points (call them $a,b$ say).
In other words blow-up $a,b,c$ and then contract the three $(-1)$-curves which are the strict transforms of the lines through $ab, bc, ac$.
The cubic becomes the desired conic; it still has 5 of the original points on it---the sixth point corresponds to the third intersection point of the line through $ab$ with the cubic. Finally the line touching the conic is the image of the exceptional curve over the cusp $c$, and this obtains two points as the images of the strict transforms of the lines through $bc$ and $ac$---the points correspond to the directions these lines went in to the cusp.

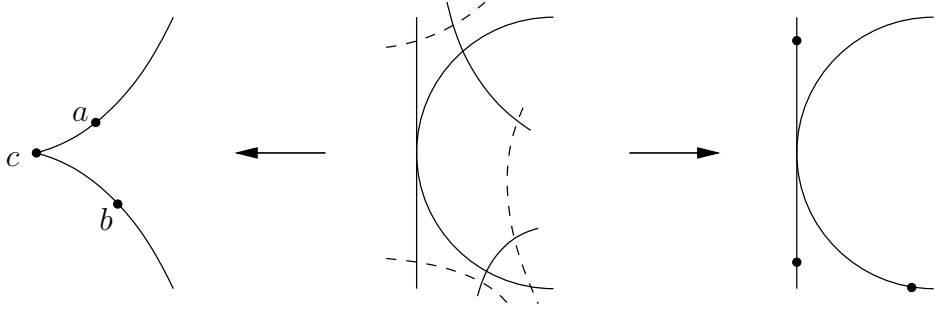
\begin{figure}[h]
\epsfxsize=5.0cm
\input{E7blowup.pstex_t}
\end{figure}

(Notice that, said differently, we have explained how to identify two eight point blow-ups of $\IP^2$, the points being either the cusp and seven points on a cuspidal cubic, or six points on a conic and two points on a line touching the conic.)

For the $E_6$ case the key point is to see how 6 points on the smooth locus of a cuspidal cubic lead to 8 points on three lines intersecting at a point, with 2 or 3 points  on each line. To do this we blow-up the cusp twice (as mentioned above to get the desired three lines) and the first three of the six points. Then we contract the following five $(-1)$-curves:

1) the strict transform of the line tangent to the cusp,

2) the strict transform of the (unique) conic through the first three points and tangent to the cubic at the cusp,

3) the strict transforms of the three lines from each of the first three points to the cusp.

Then one finds the desired 8 points on the resulting 3 lines as follows:
On the first line we still have the last 3 of the 6 original points. On the second line are the images of the three lines in step 3), and on the third line there are the images of the two $(-1)$-curves in 1) and 2).
Thus to get to Sakai's picture we need only choose one further point on the third line. (Again this can be rephrased as identifying two eight point blow-ups of $\IP^2$, and also, modulo automorphisms, there really is no further choice involved.)

\ 

\noindent{\bf Weyl group actions.\ }
Finally we will describe the Weyl group actions that occur from the present point of view (it is in these terms---i.e. Cremona transforms---that Sakai originally defined the difference equations).
Basically we will recall (from e.g. \cite{Pinkham}) 
the natural Weyl group action on $r$-tuples $p=(u_1,\ldots,u_r)$ of points of the smooth locus $U\cong \IC$ of a cuspidal cubic $C$ in $\IP^2$, for $r=6,7,8,9$. We will see how this relates to the standard Weyl group action defined in terms of the root systems and thus to the action of the previous section. (Recall that one has a canonical isomorphism $U\cong\Pic^1(C);u\mapsto [u]$ and if $0$ is the inflection point of $C$ then 
$\Pic^1(C)\cong\Pic^0(C); [u]\mapsto [u]-[0]$. In turn $\Pic^0(C)\cong \IC$ as an abelian group, although this isomorphism is only determined up to the action of 
$\IC^*$.)

Given $p=(u_1,\ldots,u_r)$ let $X=X(p)$ denote the blow-up of $\IP^2$ at the $r$ points specified by $p$.
(This makes sense even if some points coincide: one blows up the points in sequence, and at each step the remaining points determine unique points on the strict transform of the cubic, cf. \cite{Demazure-DP}.)
Then $\Pic(X)\cong \bigoplus_0^r\IZ\Eps_i$ has a uniquely determined basis
$\{\Eps_i\}$, where $\Eps_0$ is the class of the inverse image of a line in $\IP^2$, and $\Eps_i$ is the class of the {\em total} 
inverse image of $u_i\in U$. 
With respect to the intersection form on $\Pic(X)$ one has
$$\Eps_0^2=1, \quad \Eps_i^2=-1,\quad \Eps_i\cdot\Eps_j=0\quad 
\text{for all $i>0$ and  $j\ne i$.}$$
(Thus over $U^r$ we have a canonical trivialisation of the relative Picard lattice.)
The anticanonical class of $X$ is $\delta:=-\omega=3\Eps_0-\sum_1^r\Eps_i$ and is represented by the strict transform of the the cubic $C$.
Let $Q=\delta^\perp$ be the orthogonal complement of $\delta$ in $\Pic(X)$;
it is a rank $r$ lattice isomorphic to the root lattice of type
 $E_6,E_7,E_8,E_9=E_8^{(1)}$ for $r=6,7,8,9$ respectively. 
A $\IZ$-basis $\alpha$ of $Q$ is given as follows (these are the simple roots of the corresponding root system):

$$\al=\{\Eps_0-\Eps_1-\Eps_2-\Eps_3,\quad \Eps_i-\Eps_{i+1}\quad
\text{for $i=1,\ldots,r-1$}\}$$
The corresponding Weyl group acts on $\Pic(X)$ preserving $Q$, and is generated by the reflections corresponding to these roots:
$$\cF\in\Pic(X)\mapsto s_{\al_i}(\cF):= \cF + (\cF\cdot\al_i)\al_i.$$
Below, a second basis $$\be=\{\Eps_i-\Eps_0/3,\quad i=1,\ldots,r\}$$ 
of $Q\otimes \IC$ will also be useful. 

Thus, given the choice of $p$, we obtain a homomorphism 
$$\chi_p : Q\longrightarrow \Pic^0(C);\quad \cL\mapsto \cL\vert_C$$
which may be written explicitly in the $\al$ basis as
$$\chi(\Eps_0-\Eps_1-\Eps_2-\Eps_3)=3[0]-[u_1]-[u_2]-[u_3];\quad 
\chi(\Eps_i-\Eps_{i+1})=[u_i]-[u_{i+1}].$$
In other words we have a map
$$U^r\longrightarrow \Hom(Q,\Pic^0(C));\qquad p\mapsto \chi_p.$$
The key point is that this is an isomorphism and it factorises through $\Pic^0(C)^r$. Specifically we identify  $U^r$ with $\Pic^0(C)^r$ using the above isomorphism $U\cong\Pic^0(C)$ and identify $\Pic^0(C)^r$ with $\Hom(Q,\Pic^0(C))=\Hom(Q\otimes\IC,\Pic^0(C))$  using the basis $\beta$. This just amounts to observing that (the linear extension to $Q\otimes\IC$ of) $\chi_p$ satisfies
$$\chi_p(\Eps_i-\Eps_0/3) = [u_i]-[0],$$
thus motivating the definition of the basis $\be$.
In particular the Weyl group action on $Q$ naturally induces an action on 
$\Hom(Q,\Pic^0(C))$ and thus gives an action on $U^r$ via the above isomorphism.
Note that for $r=6,7,8$ the reflection hyperplanes in $U^r$ correspond to configurations $p$ such that $X(p)$ has extra $(-2)$-curves (which contract to singularities upon taking the anti-canonical model). In terms of the group law on $U\cong\IC$ these occur if 

1) two points are equal: $u_i=u_j$

2) three  points are collinear: $u_i+u_j+u_k=0$

3) six points are on a conic: $\sum_{j=1}^6 u_{i_j}=0$

4) eight points are on a cubic with a double point:  
$2 u_{i_1}+\sum_{j=2}^7 u_{i_j}=0$.

(For $r=9$ one can add any integer multiple of $\sum_1^9u_i$ to the equations for these hyperplanes---noting that this corresponds to adding multiples of 
$\delta$ to the equations for the hyperplanes in $\Hom(Q,\IC)$.)
Geometrically the generators of these Weyl group actions arise as follows. For roots $\Eps_i-\Eps_j$ one is simply swapping the labels of the points $u_i$ and $u_j$ thereby changing the order they are blown-up in. For the root 
$\Eps_0-\Eps_1-\Eps_2-\Eps_3$ one is performing the quadratic Cremona transform based at the three points $u_1,u_2,u_3$. One may verify this gives the expected change of basis of $\Pic(X)$ and that the action on $U^r$ is given (in terms of the group law on $U\cong\IC$) by
$$u_i\mapsto u_i - 2 \eps\qquad i=1,2,3,$$
$$u_i\mapsto u_i + \eps\qquad i>3$$
where $\eps=(u_1+u_2+u_3)/3$. (This may be computed explicitly---after the Cremona transform the cubic $C$ becomes another cuspidal cubic which may be mapped back to $C$ via an automorphism of $\IP^2$ and then scaled appropriately, cf.
\cite{Pinkham} pp.196-197, or more recently \cite{KMNOY} equation (139). The point is that these formulae then coincide with the reflection of $U^r$ specified above using the $\beta$ basis.)

Thus from this viewpoint it is the reflection corresponding to the single node on the shortest leg of the Dynkin diagram which is exceptional, whereas in the quiver framework it is the reflection corresponding to the central node which is unusual, as we will see in the appendix.

\begin{rmk} \label{rmrk}
Note that a priori we have defined {\em two} partial compactifications of the semi-universal deformation, one by including connections on nontrivial holomorphic bundles and one by blowing up a ninth point. One may argue they are the same, at least off of 
the affine root hyperplanes, as follows. (This was suggested by D. Arinkin.)
In both cases one has a biregular affine Weyl group action on a family of surfaces over the complement of the hyperplanes. One then checks in both cases that any point of the partial compactification divisor may be moved into the interior using this action. In terms of blowing up points this is elementary and in terms of connections this is essentially one of the statements of the Riemann--Hilbert problem actually established by Plemelj (see \cite{AnoBol94}), since off the hyperplanes at least one of the local monodromies will be semisimple. 
\end{rmk}

\end{section}

\appendix
\begin{section}{Derivation of the reflections}

The aim of this appendix is to write down how we think geometrically (in terms of connections/Fuchsian systems) of the reflections generating the finite Weyl group actions (the translation subgroup of the affine Weyl group coming from Schlesinger transformations).
In brief the reflections corresponding to the nodes on the legs of the (star-shaped affine Dynkin) quiver amount to permuting the eigenvalues of the residues of the corresponding Fuchsian system (after one normalises each residue of the system to be traceless), so for example in the $E_8$ case one immediately obtains an action of $\sym_2\times\sym_3\times \sym_6$. 
The trickier part is to obtain the reflection corresponding to the central node. In the $D_4$ (Painlev\'e VI) case this was obtained in \cite{k2p} by passing to a system of higher rank (increased by one), permuting the eigenvalues of one of the residues of this ``incremented" system and then mapping back down to the 
usual rank two framework (see also \cite{srops}).
Here we will describe the natural extension of this to the 
$E_6, E_7$ and $E_8$ cases (one still passes to systems with incremented rank, permutes the eigenvalues there and maps back down).\footnote{This is related to Katz's middle-convolution functor \cite{Katz-rls}---however the direct approach here is different to that of \cite{DettReit-katz} and should be viewed as the complex analytic analogue of Katz's approach using the l-adic Fourier transform. (We only learnt of a relation to middle convolution via \cite{Dett-Reit-p6} after writing 
\cite{k2p}---cf. diagram 1 of \cite{nlsl} for more on the origins of our approach. In turn we hope this sheds more light on Katz's functor.) Note also that \cite{CB-ihes} Theorem 3.2 independently relates the central reflection to middle convolution.} 
To deal with all the cases at the same time we will work in a more general context as follows. (In particular this gives many more examples of unexpected isomorphisms between spaces of Fuchsian systems.)

\ 

\noindent{\bf Incremented almost affine quivers.}
Let $\cQ$ be a star-shaped quiver with $m\ge 3$ legs and dimension vector $\delta$. Define an integer $N$ such that the dimension associated to the central node is $N-1$.
We will say that $(\cQ,\delta)$ is {\em almost affine} if it satisfies the following three conditions:

$\bullet$ The dimensions strictly decrease going down each leg,

$\bullet$ One of the legs is {\em full}, i.e. has dimensions $N-1,N-2\ldots,2,1$. Label the legs such that this is the last ($m$th) leg.

$\bullet$ If $n_i$ is the dimension at the node adjacent to the central node ($i=1,\ldots,m$, so $n_m=N-2$) then

\beq \label{eq: aa sum condn}
\sum_1^m n_i = 2(N-1).
\eeq

Since they are McKay graphs, any star-shaped affine Dynkin diagram (with the previously given dimension vector) satisfies \eqref{eq: aa sum condn}
and one easily observes they are almost affine.
Note also that the condition  \eqref{eq: aa sum condn} implies that the reflection $s_1$ corresponding to the central node of $\cQ$ fixes the dimension vector $\delta$.

Any such almost affine quiver may be ``incremented'' to give a new quiver as follows:
lengthen the $m$th leg by one node and increase by $1$ all the dimensions on this leg (including the dimension of the central node, which will now have dimension $N$).
Let $\cQ^+$ denote the incremented quiver obtained from an almost affine quiver $\cQ$ in this way, and let $\Delta$ be its dimension vector.

For example in the $D_4$ and $E_7$ cases the incremented quivers so obtained are:

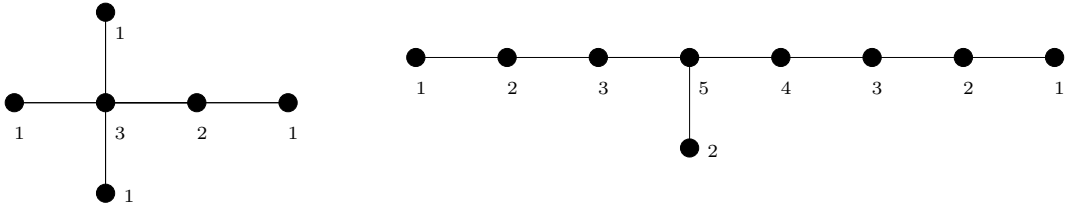
\begin{figure}[h]
\epsfxsize=5.0cm
\input{inc.dynkin.D4.E7.pstex_t}
\end{figure}

Note that an incremented quiver again has a full leg (with dimensions
$N,N-1\ldots,2,1$). Also clearly the dimension of the node adjacent to the central node on the $i$th leg 
again has dimension $n_i$, for $i=1,\ldots,m-1$. The relation 
\eqref{eq: aa sum condn} immediately implies that the central dimension $N$ partitions into these dimensions $n_i$, i.e. we have the relation:
\beq \label{eq: sum condn}
\sum_1^{m-1} n_i = N
\eeq
which will be important below.

\ 

\noindent{\bf Dual description.}\ 
Now we will give another description of the quiver varieties associated to an incremented quiver $\cQ^+$, which may be viewed as `dual' to the Fuchsian description of section \ref{sn: quiver+matrices}. 
This is a case of the duality studied by Harnad \cite{Harn94}.

Choose a set of parameters $\la$ for $\cQ^+$ (a complex scalar for each node) satisfying the condition
$$\la\cdot\Delta=0$$
where $\Delta$  is the given dimension vector of $\cQ^+$.
Then we may construct the quiver variety $N_{\cQ^+}(\la)$ as described in 
section \ref{sect, ALE}. 
Firstly, as in section \ref{sn: quiver+matrices}, $N_{\cQ^+}(\la)$
may be described as a space of Fuchsian systems; it is isomorphic to the symplectic quotient of a product of complex (co)adjoint orbits:

$$ N_{\cQ^+}(\la) \cong \wh \cO_1\times \cdots\times \wh \cO_m \spq 
\GL_{N}(\IC)$$
$$= \{ (B_1,\ldots,B_m)\ \bigl\vert\ B_i\in\wh \cO_i, \sum_1^m B_i=0 \}/\GL_{N}(\IC).$$

Here for $i<m$ we take $\wh \cO_i\subset \gl_{N}(\IC)$ to be the orbit associated to the $i$th leg of $\cQ^+$ (see section 
\ref{sn: quiver+matrices}).
For $i=m$ if we set $\cO'_m\subset \gl_{N}(\IC)$ to be the orbit associated to the last leg of $\cQ^+$, then we define  
$\wh \cO_m\subset \gl_{N}(\IC)$ to be $\cO'_m$ shifted by the scalar $\nu$ ($=$ minus the component of $\la$ associated to the central vertex of $\cQ^+$):
$$\wh \cO_m = \{ A + \nu \ \bigl\vert\ A\in \cO'_m \}.$$
Then the relation $\la\cdot \Delta=0$ translates into the statement
$\sum_1^m \tr B_i = 0$.
Note for example that the quiver varieties for the incremented $D_4$ quiver thus parameterise the $3\times 3$ Fuchsian systems with four poles used in \cite{pecr, k2p}, having rank one residues at three of the poles.
For simplicity we will impose the genericity assumption on $\la$ that 
all the orbits $\wh \cO_i$ are semisimple (and so are closed). In particular 
$\wh \cO_m$ is then regular semisimple.

For the dual description first define orbits 
$$\breve \cO_i\subset \gl_{n_i}(\IC)$$
for $i=1,\ldots,m-1$ as follows:
Remove the central vertex of $\cQ^+$ and let $\cO'_i$ be the orbit associated to what remains of the $i$th leg (the top end of which now has dimension $n_i$, so that $\cO'_i\subset \gl_{n_i}(\IC)$).
These orbits will all have zero as an eigenvalue. Define $\breve \cO_i$ by shifting by the scalar $\la_i$ (where $\la_i$ is the scalar associated to  the vertex with dimension $n_i$):
$$\breve \cO_i:=\cO'_i-\la_i\subset \gl_{n_i}(\IC).$$
(Thus for example if the full $i$th leg looked as in exercise \ref{eg: legex} of section \ref{sn: quiver+matrices} then $n_i=3$, elements of $\cO'_i$ will have eigenvalues
$0, -\la_2, -\la_2-\la_1$ and 
elements of $\breve \cO_i$ will have eigenvalues
$-\la_3, -\la_3-\la_2, \la_3-\la_2-\la_1$.)
In the $D_4$ case these orbits are single points in $\gl_1(\IC)=\IC$.

Now observe that since $\sum_1^{m-1}n_i=N$, if we define the group
$$H:=\prod_1^{m-1}\GL_{n_i}(\IC)$$ then $H$ embeds block diagonally into $\GL_N(\IC)$.
Taking Lie algebras gives the block-diagonal embedding
$$\Lia(H)=\gl_{n_1}(\IC)\oplus\cdots\oplus\gl_{n_{m-1}}(\IC)\hookrightarrow
\gl_N(\IC).$$
In particular this embeds the product of the orbits $\breve \cO_i$ into 
$\gl_N(\IC)$.
Dually projecting onto the block-diagonal entries gives a map
$\pi:\gl_N(\IC)\to \Lia(H).$
The restriction $\mu_H:=\pi\vert_{\wh \cO}$ of $\pi$ to an orbit $\wh \cO\subset \gl_N(\IC)$ is a moment map for the conjugation action of $H\subset \GL_N(\IC)$ on $\wh \cO$ (where as usual for any $k$ we identify 
$\gl_k(\IC)\cong\gl_k^*(\IC)$ using the trace pairing $\tr(AB)$).
Here we will take
$$\wh \cO := -\wh \cO_m = \{B\ \bigl\vert\ -B\in \wh \cO_m\} \subset \gl_N(\IC)$$
to be minus the regular semisimple  coadjoint orbit $\wh \cO_m$ defined above.

\begin{prop}
The quiver variety $N_{\cQ^+}(\la)$ is also isomorphic to the symplectic quotient of $\wh \cO\subset \gl_N(\IC)$ by the action of $H=\prod_1^{m-1}\GL_{n_i}(\IC)$ at the coadjoint orbit $\breve \cO:= \prod_1^{m-1}\breve \cO_i$ of $H$:
$$N_{\cQ^+}(\la) \cong \wh \cO \underset{\breve \cO}{\spq}H =
 \mu_H^{-1}\bigl( \breve \cO\bigr)/H.$$
\end{prop}
\pf
This arises since $\wh \cO \spq_{\breve \cO}H$ and 
$\prod_1^m\wh \cO_i\spq \GL_N(\IC)$ are just two descriptions of the symplectic quotient 
of $T^*\gl_N(\IC)$ by the group $H\times \GL_N(\IC)$ 
(at the coadjoint orbit $(-\breve \cO)\times \wh \cO$)---one may first do the $H$ reduction and then the $\GL_N(\IC)$ reduction to obtain one description or do the reductions  in the other order  for the other description).
Note that  $\mu_H^{-1}\bigl( \breve \cO\bigr)\subset \wh \cO$ will be empty unless the trace of an element of $\wh \cO$ equals the trace of any element of $\breve \cO$---however the equality of these two traces is equivalent to the  condition $\la\cdot \Delta=0$.

This also fits directly into the original definition of the quiver variety, as follows. Let $\cQ'\subset \cQ^+$ be the subquiver consisting of the central node (with dimension $N$) plus the $m-1$ adjacent nodes on the first $m-1$ legs (with dimensions $n_1,\ldots,n_{m-1}$ resp.), so that $\cQ'$ has exactly $m$ nodes.
Notice that the group of automorphisms $\wt G(\cQ')$ of the vector spaces at the nodes of $\cQ'$ is just $H\times \GL_N(\IC)$. Moreover
the big vector space $W(\cQ')$ associated to $\cQ'$ (taking a map in each direction for each edge): 
$$W(\cQ') = 
\bigoplus_1^{m-1} \left(\Hom(\IC^{n_i},\IC^N)\oplus\Hom(\IC^N,\IC^{n_i})\right)\cong\gl_N(\IC)\times\gl_N(\IC)\cong T^*\gl_N(\IC)
$$
is isomorphic to $ T^*\gl_N(\IC)$, once we make an 
identification $\bigoplus_1^{m-1} \IC^{n_i}\cong \IC^N$ (as we did above embedding $H$ in $\GL_N(\IC)$ block diagonally).

Now we may do the symplectic quotient 
$N_{\cQ^+}(\la)=W(\cQ^+)\spq_\la\wt G(\cQ^+)$ in stages, first quotienting by the subgroup $K$ of $\wt G(\cQ^+)$ corresponding to the nodes of $\cQ^+\setminus \cQ'$, so that $\wt G(\cQ^+)=K\times H\times\GL_N(\IC)$.
This yields the isomorphism:
$$W(\cQ^+)\spq_\la K \cong \breve \cO_1\times\cdots\times\breve \cO_{m-1}\times \wh \cO_m\times T^*\gl_N(\IC)=\breve \cO\times \wh \cO_m\times T^*\gl_N(\IC)$$ 
on which we still have the residual action of $H\times\GL_N(\IC)$.
Suppose $P_i\in \Hom(\IC^N,\IC^{n_i})$ and $Q_i\in \Hom(\IC^{n_i},\IC^N)$ for $i=1,\ldots,m-1$ and let us write
$$P:= 
\left(\begin{matrix} 
P_1 \\
\vdots \\
P_{m-1}
\end{matrix}\right), \quad Q:= (Q_1,\ldots,Q_{m-1})
\in \gl_N(\IC).
$$
Then $(P,Q)\in T^*\gl_N(\IC)$ and $\GL_N(\IC)$ acts on this via left or right multiplication on $\gl_N(\IC)$. Specifically $g\in \GL_N(\IC)$ acts as
$g(P,Q) = (Pg^{-1}, gQ)$
with moment map
$$\mu_L(P,Q)= QP=\sum_1^{m-1} Q_iP_i$$
and $h\in H\subset \GL_N(\IC)$ acts as 
$h(P,Q) = (hP, Qh^{-1})$ with moment map
$$\mu_R(P,Q)= -\pi(PQ) = (-P_1Q_1,-P_{2}Q_{2},\ldots)\in\gl_{n_1}(\IC)\times\gl_{n_2}(\IC)\times\cdots.$$
Thus we have immediately that
$$N_{\cQ^+}(\la)\cong 
(\breve \cO\times \wh \cO_m\times T^*\gl_N(\IC))\spq (H\times\GL_N(\IC)) $$
$$\cong\{ (P,Q)\in T^*\gl_N(\IC)\ \bigl\vert\ 
P_iQ_i\in \breve \cO_i, QP \in \wh \cO \}/(H\times \GL_N(\IC))$$
$$\cong \mu_R^{-1}(-\breve \cO)\cap \mu_L^{-1}(\wh \cO)/(H\times \GL_N(\IC))$$
\beq\label{eq: last im}
=T^*\gl_N(\IC) \underset{(-\breve \cO)\times\wh \cO}\spq (H\times GL_N(\IC)).
\eeq

\vspace{2mm}
\noindent
Now we claim that there is an isomorphism 
$$T^*\gl_N(\IC) \spq_{\wh \cO}\GL_N(\IC)\cong \wh \cO.$$
(which will be an isomorphism of symplectic manifolds upon negating the symplectic form on the left-hand side).
This claim will complete the proof of the proposition since, from 
\eqref{eq: last im},
we then obtain 
$N_{\cQ^+}(\la)\cong \wh \cO\spq_{\breve \cO} H$.

To establish the claim we need to identify the affine GIT quotient
$$\{ (P,Q)\in T^*\gl_N(\IC)\ \bigl\vert\ QP\in\wh \cO\}/\GL_N(\IC)$$
with $\wh \cO$. 
The first step is to note that the invariant functions are the matrix entries of $B:=PQ\in\gl_N(\IC)$.
Now if zero is not an eigenvalue of elements of $\wh \cO$, then $\GL_N(\IC)$ acts freely and $PQ$ is conjugate to $QP$, and the result follows easily.

Next suppose zero is an eigenvalue of elements of $\wh \cO$, and we have $(P,Q)$ with $QP\in\wh \cO$, so $\det(QP)=\det(PQ)=0$.
Label the nonzero eigenvalues of $QP$ as $\nu_1,\ldots,\nu_{N-1}$ and let $E_i\subset \IC^N$ be the $\nu_i$-eigenspace of $QP$ for $i=1,\ldots,N-1$, and let $E_0=\ker(QP)$
(thus $E_i\cong\IC$ as $\wh \cO$ is regular semisimple).
Let $F_i:= P(E_i)\subset \IC^N$ for $i>0$ and observe $F_i\cong \IC$, $F_i$ is the $\nu_i$-eigenspace of $PQ$, and that $Q(F_i)=E_i$.
Let $F_0:=\ker(PQ)$. Counting dimensions implies $\dim(F_0)=1$ and so 
$PQ\in\wh \cO$.

Thus the map $(P,Q)\mapsto B:=PQ$ is a $\GL_N(\IC)$-invariant map to $\wh \cO$, which is clearly surjective.
It is now straightforward to show that each fibre of this map contains exactly three $\GL_N(\IC)$ orbits:
two non-closed orbits (with $\det(P)=0\ne\det(Q)$ and $\det(P)\ne 0=\det(Q)$) and one closed orbit (with $\det(P)=\det(Q)=0$) in the closure of both of the non-closed orbits.
Thus the GIT quotient (identifying points if the closures of their orbits intersect) is $\wh \cO$.
\epf

\begin{rmk}\label{rmk: closed orbit}
For use later we note that if $(P,Q)\in T^*\gl_N(\IC)$ and $QP\in\wh \cO$ (and $\det(QP)=0$) then we have shown that the closure of the $\GL_N(\IC)$ orbit of $(P,Q)$ contains elements with $\det(P)=\det(Q)=0$---and so any point of the GIT quotient may be represented by such $(P,Q)$.
Now if $\det(P)=\det(Q)=0$ we may deduce that
$$\ker(P)=\ker(QP)=:E_0,\qquad \ker(Q)=\ker(PQ)=:F_0.$$
(E.g. for the first statement, a dimension count shows the image of $P$ is $\bigoplus_1^{N-1}F_i$. Then if $v\in E_0$ we have $P(v)=\sum_1^{N-1}w_i$ with $w_i\in F_i$. Thus $0=QPv=\sum_1^{N-1}v_i$ with $v_i=Qw_i\in E_i$, and $v_i=0$ if and only if $w_i=0$. But the $E_i$'s sum directly so each $v_i=0$.)
It then follows that, if we define $V=\bigoplus_1^{N-1}E_i\cong\IC^{N-1}$, then $V\oplus E_0=\IC^N$ and each matrix $B_i:=Q_iP_i$ preserves this direct sum decomposition (since $E_0\subset \ker P_i$ and the image of $Q_i$ is contained in the image $Q(\IC^N)=V$ for all $i$).
\end{rmk}

\begin{rmk}
As in \cite{Harn94, smapg} we may interpret 
$\wh \cO\spq_{\breve \cO} H$ 
as a moduli space of certain irregular connections of the form
$$d-\left(\frac{A_0}{z^2} + \frac{B}{z}\right)dz$$
on the trivial rank $N$ vector bundle over $\IP^1(\IC)$, where 
$B=PQ$ is in $\wh \cO$ and $A_0\in \gl_n(\IC)$ is a diagonal element with eigenvalues of multiplicities $n_1, n_2,\ldots, n_{m-1}$ respectively (and so having centraliser $H$).
\end{rmk}

\ 

\noindent{\bf Scalar shift.}\ 
Now we consider the effect of tensoring the above irregular connections by a logarithmic connection $\Lambda dz/z$ on the trivial line bundle over $\IP^1$. In other words we replace $B$ by $B+\Lambda$ where $\Lambda\in\IC$ is a scalar.
The effect on the parameters $\la$ is clear since we see how the eigenvalues of $B$ change: only the components of $\lambda$ associated to the nodes of $\cQ'$ change, viz $\lambda_0$ (associated to the central node with dimension $N$) increases by $\Lambda$ and each $\lambda_i$ (associated to the node with dimension $n_i$) decreases by $\Lambda$. Let us denote this new set of parameters by $\la(\Lambda)$, so that $\la=\la(0)$. Note that the relation $\sum n_i=N$ implies we have $\la(\Lambda)\cdot\Delta=0$ for all $\Lambda$.
The key point is that this map $B\mapsto B+\Lambda$ induces an isomorphism
$$N_{\cQ^+}(\la)\cong N_{\cQ^+}(\la(\Lambda))$$
for any $\Lambda\in\IC$. This is immediate in this dual description but is less trivial in the original Fuchsian description---indeed for generic $\Lambda$ this is (our approach to) Katz's middle convolution. (At least if zero is not an eigenvalue of elements of $\wh \cO$, one may describe this scalar shift on the Fuchsian side by using the $\GL_N(\IC)$ action to move to points where $Q=\id$---then we see how to vary $P=B$ with $\Lambda$ and in turn how $B_i=Q_iP_i\in \wh \cO_i(\Lambda)$ varies. In general one may move to a GIT-equivalent pair with $Q=1$, and proceed in the same way.)

\ 

\noindent{\bf Projection to original quiver variety.}\ 
Suppose $(\cQ,\delta)$ is an almost affine quiver with incremented quiver $(\cQ^+,\Delta)$ and that we have chosen some 
generic parameters $\la$ for $\cQ^+$.
Let $\nu$ be minus the component of $\lambda$ associated to the central vertex of $\cQ^+$.
Then we may shift $\la$ by $\Lambda=\nu$ so that the central component $\la(\Lambda)$ associated to the central node of $\cQ^+$ is zero. Thus $\la(\Lambda)$ constitutes a set of parameters for the original quiver $\cQ$ (since $\cQ$ may be obtained by deleting the central vertex of $\cQ^+$ and then contracting the highest joint of the $m$th leg of $\cQ^+$). Writing $\pr(\la)$ for the set of parameters of $\cQ$ obtained from $\la$ in this way, we have
$$\pr(\la)\cdot\delta = \la(\Lambda)\cdot\Delta = \la\cdot\Delta = 0.$$

\begin{lem}
The quiver variety of $\cQ^+$ with parameters $\la$ is isomorphic to the quiver variety of $\cQ$ with parameters $\pr(\la)$:
$$N_{\cQ^+}(\la)\cong N_\cQ(\pr(\la)).$$
\end{lem}
\sketch
By performing the scalar shift we may assume that 
the component of $\la$ on the central node is zero.
Then if $B=PQ\in \wh \cO$ represents a point of $N_{\cQ^+}(\la)$ we know that $B$ has exactly one zero eigenvalue. 
From Remark \ref{rmk: closed orbit} we can move to a GIT-equivalent pair $(P,Q)$ with $\det(P)=\det(Q)=0$ and then the matrices $B_i=Q_iP_i$ simultaneously preserve a direct sum decomposition of the form 
$\IC^N=V\oplus\IC$. Taking $A_i=B_i\vert_V\in \cO_i$ gives the desired set of 
Fuchsian residues representing a point of $N_\cQ(\pr(\la))$.
(In reverse we may take $B_i\in\wh \cO_i$ to be the block diagonal matrix 
$\diag(A_i,0)\in \gl_N(\IC)$.)

\esketch

\begin{eg}
To help keep track of the various orbits we have defined let us write out the eigenvalues in the example of $E_7$. If we number the nodes of $\cQ^+$ as follows
\begin{figure}[h]
\epsfxsize=5.0cm
\input{inc.dynkin.E7.labels.pstex_t}
\caption{}\label{fig: inc e7 labels}
\end{figure}

\noindent
so that the first leg is the shortest leg, and the longest leg is last, then
the eigenvalues of the various orbits are:

\noindent$\bullet$
Any element of $\wh \cO_1$ has eigenvalues $(0,0,0,-\la_1,-\la_1)$
and
$\breve \cO_1 = \{ -\diag(\la_1,\la_1) \}\subset \gl_2(\IC)$,

\noindent$\bullet$
Any element of $\wh \cO_2$ has eigenvalues 
$(0,0,-\la_2,-\la_2-\la_7, -\la_2-\la_7-\la_8)$ and
any element of $\breve \cO_2\subset\gl_3(\IC)$ has eigenvalues $(-\la_2,-\la_2-\la_7, -\la_2-\la_7-\la_8)$.

\noindent$\bullet$
Any element of $\wh \cO_3$ has eigenvalues 
$-\la_0$ and $-\la_0-\sum_{i=3}^j\la_i$ for $j=3,4,5,6$
and any element of $\wh \cO$ is minus an element of $\wh \cO_3$.

The parameters of the original affine Dynkin diagram obtained by projecting $\la$ are then as follows:
\begin{figure}[h]
\epsfxsize=5.0cm
\input{dynkin.E7.labels.pstex_t}
\caption{}\label{fig: e7 projn}
\end{figure}
\end{eg}

Thus we have identified $N_\cQ(\pr(\la))$ with $N_{\cQ^+}(\la(\La))$ for any $\La\in\IC$, and so with spaces of Fuchsian systems of the form
$$\wh \cO_1\times\cdots\times\wh \cO_m \spq \GL_N(\IC)$$
of incremented rank.
In particular the size of the `obvious' symmetry group (permuting the eigenvalues of the residues) is larger---namely we may permute all $N$ of the eigenvalues of elements of $\wh \cO_m$ compared to the $N-1$ eigenvalues of $\cO_m$. In particular suppose we consider the permutation swapping the first two eigenvalues of elements of $\wh \cO_m$. If the parameters on the $m$th leg of $\cQ^+$ were labelled  as follows:
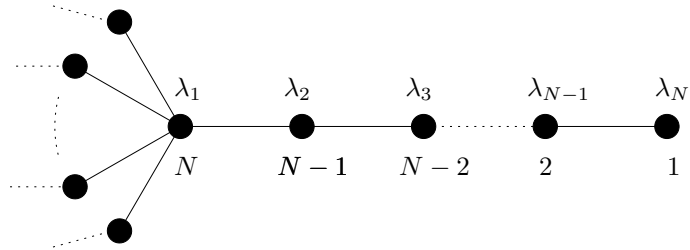
\begin{figure}[h]
\epsfxsize=5.0cm
\input{quiverleg.pstex_t}
\end{figure}

\noindent
then the eigenvalues would be $-\sum_{j=1}^i\la_j$ for $i=1,2,\ldots,N$.
Let $\per(\la)$ denote the change in the parameters $\la$ corresponding to this permutation 
swapping the first two eigenvalues. Explicitly we find that
$$\per(\la):\qquad\la'_1=\la_1+\la_2,\quad \la'_2=-\la_2,\quad 
\la'_3= \la_2+\la_3$$
with all other components of $\la$ not changing.
The key fact is that this permutation is a lift of the desired `central reflection' of $\cQ$, as follows.

\begin{lem}
Let $r_1$ denote the reflection corresponding to the central node of the almost affine quiver $\cQ$.
Then, upon projecting to $\cQ$, the eigenvalue permutation $\per$ induces the reflection $r_1$. In other words:
$$r_1(\pr(\la)) = \pr(\per(\la))$$
for any set of parameters $\la$ for $\cQ^+$ (with $\la\cdot \Delta=0$).
\end{lem}
\pf
Let $\nu$ be minus the component of $\pr(\la)$ associated to the central 
vertex of $\cQ$. 
Then $r_1$ acts by subtracting $\nu$ from all the components of $\pr(\la)$  on nodes of $\cQ$ adjacent to the central node (and negating $\nu$ itself).
(In terms of eigenvalues---in the `determinant zero' normalisation---this amounts to subtracting $\nu$ from all of the not-necessarily-zero eigenvalues of each residue. E.g. for $D_4$ one has $\nu=\sum_1^4 \theta_i/2$ and $r_1$ maps each $\theta_i$ to $\theta_i-\nu$.)

The proof is now straightforward---we will illustrate it for the $E_7$ case, the general case being no more difficult.

From Figure \ref{fig: e7 projn} we see that $r_1(\pr(\la))$ is
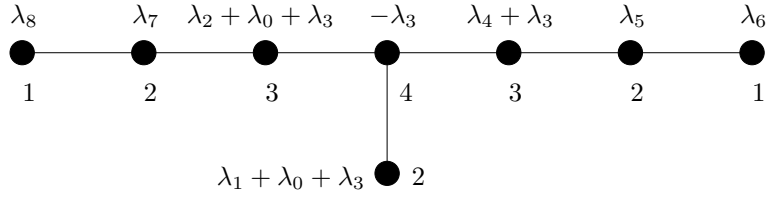
\begin{figure}[h]
\epsfxsize=5.0cm
\input{dynkin.E7.r1labels.pstex_t}
\end{figure}

Moreover from Figure \ref{fig: inc e7 labels}, we see that $\per(\la)$ is
as follows:

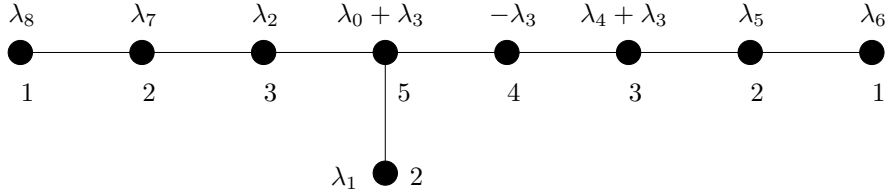
\begin{figure}[h]
\epsfxsize=5.0cm
\input{inc.dynkin.E7.plabels.pstex_t}
\end{figure}

Now observe that these parameters project to those above, as desired.
\epf

Of course relabelling the eigenvalues of elements of $\wh \cO_m$ does not change the spaces, so we have an isomorphism
$$N_{\cQ^+}(\la)\cong N_{\cQ^+}(\per(\la))$$
(which is transparent when the quiver varieties are viewed as spaces of Fuchsian systems).
Putting the above results together, we thus see how to obtain the desired isomorphism $N_\cQ(\pr(\la))\cong N_\cQ(r_1(\pr(\la)))$ as the composition
$$N_\cQ(\pr(\la))\cong N_{\cQ^+}(\la)\cong N_{\cQ^+}(\per(\la))\cong
N_\cQ(\pr(\per(\la)))=
N_\cQ(r_1(\pr(\la)))$$
as promised by passing to the spaces of incremented systems, permuting eigenvalues and projecting back down.

\end{section}

\ 

\noindent
{\bf Acknowledgments.}

\noindent
The author is grateful to E. Looijenga for useful conversations and to D. Arinkin for Remark \ref{rmrk}.

\renewcommand{\baselinestretch}{0.9}              %
\normalsize
\bibliographystyle{amsplain}    \label{biby}
\bibliography{../thesis/syr}    
\end{document}